\documentclass[12pt]{article}
\usepackage{theorem,amsfonts}

\newtheorem{theorem}{Theorem}[section]
\newtheorem{proposition}{Proposition}[section]
\newtheorem{lemma}{Lemma}[section]
\newtheorem{corollary}{Corollary}[section]
\theorembodyfont{\upshape}
\newtheorem{definition}{Definition}

\newtheorem{remark}{Remark}[section]
\newtheorem{proof}{Proof}

\newcommand{\bt}{\begin{theorem}}
\newcommand{\et}{\end{theorem}}
\newcommand{\bl}{\begin{lemma}}
\newcommand{\el}{\end{lemma}}
\newcommand{\bp}{\begin{proposition}}
\newcommand{\ep}{\end{proposition}}
\newcommand{\bd}{\begin{definition}}
\newcommand{\ed}{\end{definition}}
\newcommand{\br}{\begin{remark}}
\newcommand{\er}{\end{remark}}
\newcommand{\bc}{\begin{corollary}}
\newcommand{\ec}{\end{corollary}}
\newcommand{\bo}{\begin{proof}}
\newcommand{\eo}{\end{proof}}
\newcommand{\be}{\begin{enumerate}}
\newcommand{\ee}{\end{enumerate}}

\title{No representation of Moore groups and affine groups has any rate of 
random mixing}
\author{C. R. E. Raja}
\date{ }

\begin{document}
\maketitle

\let\epsi=\epsilon
\let\lam=\lambda
\let\ra=\rightarrow
\let\da=\downarrow
\let\hra=\hookrightarrow
\let\da=\downarrow 
\let\ss=\subset
\let\ba=\beta
\let\ga=\gamma
\let\Ga=\Gamma

\begin{abstract}
A sequence $a_n\da 0$ forms a rate of random mixing for an unitary system 
$(G,\mu , \pi, {\cal H})$ if for any $u, v\in {\cal H}$ 
$$\limsup {|<\pi (g_n^\omega )u, v>|\over a_n} < \infty$$ a.e. $\omega$ in 
the probability space $(G^{\mathbb N}, \mu ^{\mathbb N})$ of the random 
walk induced by $\mu$.  We study the class of locally compact groups 
none of whose representation has any rate of random mixing and prove that 
this class contains Moore groups and certain solvable groups which 
includes the group of affine transformations on a local field.    
\end{abstract}



\begin{section}{Introduction and Preliminaries}

Let $G$ be a locally compact group and ${\cal P}(G)$ be the space of 
regular Borel probability measures on $G$ with weak topology which is 
the coarsest topology for which the functions $\mu \mapsto \mu (f)$ are 
continuous for all continuous bounded function $f$ on $G$.  

For $\mu \in {\cal P}(G)$ and $n \geq 1$, let $\mu ^n$ denote the $n$-th 
convolution product of $\mu$ with itself.  

A probability measure $\mu \in {\cal P}(G)$ is called {\it adapted} if it 
is not supported on a proper closed subgroup and $\mu$ is called {\it 
strictly aperiodic} if it is not supported on a coset of a proper closed 
normal subgroup.  

A probability measure $\mu \in {\cal P}(G)$ is called {\it spread-out} if 
$\mu ^k$ is not singular with respect to a Haar measure on $G$ for some 
$k \geq 1$.  

Let $\mu$ be a probability measure on $G$.  Let $G^{\mathbb N}$ be denote 
the product space $\Pi _{i=1}^\infty G$ which is the space of paths 
$\omega = (\omega _n)$ of the random walk and $\mu ^{\mathbb N}$ be the 
product measure $\Pi _{i=1}^\infty \mu$.  Then $(G^{\mathbb N}, \mu 
^{\mathbb N})$ is called the probability space of the random walk induced 
by $\mu$.  For any sequence $\omega = (\omega _n)\in G^{\mathbb N}$, we 
define the $n$th random product by 
$$g_n^\omega = w_n w_2 \cdots w_1$$ for any $n \geq 1$.  

We say that $(G, \mu , \pi, {\cal H})$ is an unitary system if $G$ 
is a locally compact group, $\mu \in {\cal P}(G)$ is adapted and 
strictly aperiodic spread-out probability measure and $\pi$ is an unitary 
representation of $G$ on a Hilbert space $\cal H$.  

\bd
Let $(G, \mu, \pi, {\cal H})$ be an unitary system.  A decreasing 
sequence $(a_n)$ with $a_n \ra 0$ is said to form a {\it rate of 
random mixing} for $(G, \mu, \pi, {\cal H})$ if for every 
$u, v\in {\cal H}$, with probability one: 
$$\limsup {|<\pi (g_n^\omega )u, v>| \over a_n} <\infty .$$
\ed

The notion of rate of random mixing was introduced in \cite{Sh} to 
bound the critical exponential rate of mixing from below.  
Since every irreducible representation of a locally compact abelian group 
is one-dimensional and by Theorem 2.13 of \cite{Sh}, we get that no 
unitary representation of a locally compact abelian group has any rate of 
random mixing.  Shalom raises the question of characterizing groups with 
this property.  It follows from Corollary 5.5 of \cite{Be} and Theorem 2.8 
of \cite{Sh} that there are unitary representations of non-amenable groups 
having a rate of random mixing.  

In this note we prove that no unitary representation of Moore groups (that 
is, groups for which all irreducible unitary representations are 
finite-dimensional) and certain solvable groups (which includes the group 
of affine transformations on a local field of characteristic zero) has any 
rate of random mixing.

We now prove the following useful lemmas first of which is easy to verify 
and we omit the proof.   

\bl\label{l1}
Let $(G, \mu, \pi , {\cal H})$ be a unitary system and $H$ be a normal 
subgroup of $G$ contained in the kernel of $\pi$.  Then $(a_n)$ is a rate 
of random mixing for $(G, \mu , \pi , {\cal H})$ if and only if $(a_n)$ is 
a rate of random mixing for $(G/H, \tilde \mu , \tilde \pi , {\cal H})$ 
where $\tilde \mu$ is the projection of $\mu$ onto $G/H$ and $\tilde \pi$ 
is the factor representation of $\pi$.
\el

\bl\label{l2}
Let $(G, \mu , \pi, {\cal H})$ be a unitary system with $\pi$ 
irreducible.  Then there exists a 
compact normal subgroup $K$ of $G$ such that $K$ is contained in the 
kernel of $\pi$ and $G/K$ is second countable and $(a_n)$ is a rate of 
random mixing for $(G, \mu , \pi , {\cal H})$ if and only if $(a_n)$ is a 
rate of random mixing for $(G/K, \tilde \mu , \tilde \pi , {\cal H})$ 
where $\tilde \mu$ is the projection of $\mu$ onto $G/K$ and $\tilde \pi$ 
is the factor representation of $\pi$.
\el

\bo
Since $\mu$ is adapted, $G$ is $\sigma$-compact and hence $G$ can be 
approximated by second countable groups.  Let $v \in \cal H$ be a unit 
vector.  Then there exists a compact normal subgroup $K$ of $G$ such that 
$G/K$ is second countable and $||\pi (g) v-v|| < 1$ for all $g \in K$.  
This implies that $||\pi (\omega _K )v -v|| <1$ where $\omega _K$ is the 
normalized Haar measure on $K$.  Thus, $\pi (\omega _K)v$ is a 
$K$-invariant non-zero vector in $\cal H$.  Since $K$ is normal, the space 
of $K$-invariant vectors is an invariant subspace.  Since $\pi$ is 
irreducible, $K$ is contained in the kernel of $\pi$.  Now the second part 
follows from Lemma \ref{l1}.
\eo

In view of Lemma \ref{l2}, we may assume that the class of groups under 
consideration consists of second countable groups.  

\end{section}

\begin{section}{Moore groups}

\bp\label{f} 
Let $(G, \mu , \pi , {\cal H})$ be any unitary system.  Suppose the 
unitary representation $\pi$ is of finite-dimension.  Then for $u$ and $v$ 
in $\cal H$, $<\pi (g_{k_n}^\omega ) u , v> \ra 0$ a.e $\omega$ for some 
subsequence $(k_n)$ of positive integers if and only if there exist 
orthogonal invariant subspaces $U$ and $V$ of $\cal H$ such that $u \in U$ 
and $v \in V$.
\ep

\br
We would like to remark that Proposition \ref{f} is strictly stronger than 
not admitting any rate of random mixing and it may be seen by showing that 
the only if part of Proposition 2.1 is not true for the left regular 
representation of non-compact groups.  Since the left regular 
representation of amenable groups weakly contain the trivial 
representation, left regular representation of amenable groups have no 
rate of random mixing.  Let $G$ be a non-compact locally compact group and 
$R$ be the left regular representation of $G$.  Suppose $\mu$ is any 
adapted and strictly aperiodic spread-out probability measure on $G$ and 
$f \in L^2(G)$ be any non-negative function.  Then 
$$\int <R(g_n^\omega )f, f> d\mu ^{\mathbb N} (\omega ) = 
\int <R(g) f, f> d\mu ^n (g) = <R(\mu )f, f> \ra 0$$ by Theorem 2.8 of 
\cite{DL}.  Since $<R(g)f, f>\geq 0$ for all $g \in G$, there exists a 
subsequence $(k_n)$ such that $$<R(g_{k_n}^\omega )f , f> \ra 0$$ a.e. 
$\omega $ (see for example Theorem 3.12 of \cite{Ru}).
\er

\bo 
Let $\cal U$ be the group of unitary operators on $\cal H$.  Since 
$\cal H$ is of finite-dimension, $\cal U$ is compact.  Let $K$ be the 
closure of $\pi (G)$.  Then $K\ss {\cal U}$ and hence $K$ is compact.  

Suppose there exists a subsequence $(k_n)$ of positive integers and 
vectors $u$ and $v$ in ${\cal H}$ such that 
$<\pi (g_{k_n}^\omega ) u , v> \ra 0$ a.e. $\omega$.  For 
$n \geq 1$, let $f_n \colon G^{\mathbb N} \ra {\mathbb R}$ be defined by 
$$f_n (\omega ) = |<\pi(g_{k_n}^\omega )u, v>|$$ for all $\omega \in 
G^{\mathbb N}$.  Then $f_n$ is uniformly bounded and $f_n(w) \ra 0 $ a.e. 
$\omega$.  By Lebesque dominated convergence theorem, 
$$\int |<\pi (g)u, v>| d\mu ^{k_n}(g) = \int f_n(w) d\mu ^{\mathbb N} 
(\omega )\ra 0$$ as $n \ra \infty$.  

Let $f \colon K \ra {\mathbb R}$ be defined by $$f(t) = |<t(u), v>|$$ for 
all $t \in K$.  Then $f$ is a continuous bounded function and 
$$\int f(t) d\pi (\mu )^{k_n} (t) = \int f(\pi (g)) d\mu ^{k_n}(g) \ra 0$$ 
as $n\ra \infty$.  Since $\mu$ is adapted and strictly aperiodic in $G$, 
$\pi (\mu )$ is also adapted and strictly aperiodic in $K$.  Thus, by 
Kawada-Ito theorem  $\pi (\mu )^n  \ra \omega _K$ in ${\cal P}(K)$ where 
$\omega _K$ is the normalized Haar measure on $K$ (see \cite{KI}).  Since 
$f$ is a continuous bounded function on $K$, we get that 
$$\int  f(t) d\pi (\mu )^n (t) \ra \int f(t) d\omega _K(t)$$ as 
$n \ra \infty$.  Thus, $$\int f(t) d\omega _K(t) = 0.$$  Since $f\geq 0$ 
is a continuous function on $K$, $f(t) =0$ for all $t \in K$.  This in 
particular, implies that $<\pi (g) u, v>=0$ for all $g \in G$.  Let $U$ be 
the subspace of $\cal H$ spanned by $\{ \pi (g) u \mid g\in G\}$.  Then 
$U$ is a closed invariant subspace of $\cal H$ such that $u \in U$ and 
$v \in U^\perp$.  Thus, vectors $u$ and $v$ are in two orthogonal 
invariant subspaces.  Converse is easy to prove.  
\eo

We now apply Proposition \ref{f} to finite-dimensional representations and 
in particular to Moore groups.  A locally compact group $G$ is called a 
{\it Moore group} if irreducible unitary representations of $G$ are 
finite-dimensional.

\bt\label{m}
Let $(G, \mu , \pi , {\cal H})$ be any unitary system.  If $\pi$ is 
finite-dimensional, then $\pi$ has no rate of random mixing.  In 
particular, no representation of a Moore group has any rate of random 
mixing. 
\et

\bo
Suppose the sequence $(a_n)$ forms a rate of random mixing for $(G, \mu, 
\pi , {\cal H} )$.  Since $\pi$ is finite dimensional, $\cal H$ 
contains a non-trivial irreducible subspace and let $u$ and $v$ be from a 
non-trivial irreducible subspace.  Then 
$$\limsup { |<\pi (g_n^\omega) u, v>| \over a_n} <\infty$$ for a.e. 
$\omega $.  This implies that $|<\pi (g_n^\omega )u, v>| \ra 0$ a.e. 
$\omega$.  By Proposition \ref{f}, $u$ and $v$ are from two invariant 
orthogonal subspaces.  This is a contradiction.  Hence the unitary system 
$(G, \mu , \pi , {\cal H})$ has no rate of random mixing.

Suppose $G$ is a Moore group.  Let $(G, \mu, \pi , {\cal H})$ be an 
unitary system.  Using Theorem 2.13 of \cite{Sh}, we may assume that $\pi$ 
is an irreducible unitary representation of $G$.  Then ${\cal H}$ is 
of finite-dimension.  By the previous case, we get that 
$(G, \mu, \pi, {\cal H})$ has no rate of random mixing.
\eo

\end{section}

\begin{section}{Affine groups}

We now show that no representation of affine groups has any rate of random 
mixing.

\bt\label{af}
Let ${\mathbb K}$ be a local field of characteristic zero.  
Let $G$ be the semidirect product of ${\mathbb K}^*$, the multiplicative 
group of non-zero elements in ${\mathbb K}$ and ${\mathbb K}^n$ where 
action of ${\mathbb K}^*$ on ${\mathbb K}^n$ is given by 
$$a\cdot (a_1 , \cdots ,a_n)=(aa_1, \cdots ,aa_n)$$ for all 
$a\in {\mathbb K}^*$ and $(a_1, \cdots , a_n ) \in {\mathbb K}^n$.  Then 
no representation of $G$ has any rate of random mixing.
\et

\bo
Let $\pi$ be an irreducible unitary representation of $G$.  Then by 
Mackey's Theorem there exists a character $\chi$ of ${\mathbb K}^n$ such 
that $\pi$ is induced from an irreducible representation, say $\rho$, of 
the stabilizer $G_\chi$ of $\chi$ (see for example Theorem 6.38 of 
\cite{Fo}).  If $\chi$ is trivial, then $G_\chi = G$ and hence $\pi$ is an 
one-dimensional representation.  Thus, $\pi$ has no rate of random mixing.  

If $\chi$ is non-trivial, then $G_\chi = {\mathbb K}^n$ and $\pi$ is 
induced from the one-dimensional representation $\chi$ of ${\mathbb K}^n$.  
Let $dm$ be a Haar measure on ${\mathbb K}^*$.  Then $\pi$ is the 
representation defined on $L^2({\mathbb K}^* )$ by
$$\pi (x) f(b) = \chi (b^{-1} \cdot u) f(a^{-1}b)$$ for all $f \in L^2 
({\mathbb K}^*)$, $x = (a, u) \in G$ and $b \in {\mathbb K}^*$.    

Let $V_0$ be the subspace of ${\mathbb K}^n$ of co-dimension one such that 
$\chi =1$ on $V_0$.  Now for $f \in L^2({\mathbb K}^*)$ and $v\in V_0$, 
$$\pi (v)f(b)=\chi (b^{-1}\cdot v)f(b)=f(b)$$ 
for all $b\in {\mathbb K}^*$.  This shows that $\pi$ is trivial on $V_0$.  
Replacing ${\mathbb K}^n$ by ${\mathbb K}^n /V_0$, we may 
assume that $n =1$.  Thus, the group $G$ is the semidirect product of 
${\mathbb K}^*$ and ${\mathbb K}$.  

We now claim that $\pi$ weakly contains the trivial representation of $G$.  
Since ${\mathbb K}^*$ is amenable, there exists a summing sequence of 
non-null compact sets $(B_n)$ in ${\mathbb K}^*$ (see 4.15 and 4.16 of 
\cite{Pa}).  For $n \geq 1$, define 
$$f_n = {1\over \sqrt {dm(B_n)}}1_{B_n}$$ where 
$1_{B_n}$ is the indicator function on the set $B_n$.  Then $(f_n)$ is 
asymptotically translation invariant in $L^2({\mathbb K}^*)$, that is, 
$$||\pi (a) f_n -f_n|| = ||R(a)f_n-f_n||\ra 0$$ as $n \ra \infty$ for all 
$a\in {\mathbb K}^*$ where $R$ is the left regular representation of 
${\mathbb K}^*$ on $L^2({\mathbb K}^*)$.  Now for $u \in {\mathbb K}$, 
$$\begin{array}{rcl}
||\pi (u) f_n - f_n ||^2 & = & 
{1\over dm(B_n)}\int _{B_n} |\chi (b^{-1}u) -1| ^2 dm(b) \\
& \leq & {1\over dm(B_n)} \int _{B_n \setminus B_k} |\chi (b^{-1}u) 
-1| ^2 dm (b) + 2{dm(B_k)\over dm(B_n)} \\
\end{array}$$ for any $k \leq n$.  Since $\chi (b^{-1}u) \ra 1$ as 
$b\ra \infty$, we get that $$||\pi (u) f_n -f_n ||^2 \ra 0$$ as 
$n\ra \infty$ for all $u \in {\mathbb K}$.  Thus, the group 
$\{ g \in G\mid ||\pi (g)f _n -f_n||\ra 0\}$ contains $\mathbb K$ 
and ${\mathbb K}^*$ and hence $||\pi (g)f_n-f_n||\ra 0$ for all $g\in G$.  
This implies that the trivial representation is weakly contained in $\pi$ 
(see Remark 2.1 of \cite{Ra}).  By Theorem 2.13 of \cite{Sh}, any rate of 
random mixing of $\pi$ is also a rate of random mixing of the trivial 
representation but the trivial representation has no rate of random 
mixing.  Hence $\pi$ has no rate of random mixing. 
\eo

\end{section}

\noindent{C. Robinson Edward Raja \\ Stat-Math Unit \\ Indian Statistical 
Institute \\ 8th Mile Mysore Road \\ Bangalore -560 059. \\ e-mail: 
creraja@isibang.ac.in}

\end{document}